\documentclass{bmcart}

\usepackage{amsthm,amsmath}
\usepackage[utf8]{inputenc} 
\usepackage{amssymb,bm}
\usepackage{amsfonts}
\usepackage{graphicx}
\usepackage{afterpage}
\usepackage{flafter}
\usepackage{marvosym}
\usepackage{multirow}
\usepackage{dcolumn}
\usepackage{lipsum}
\usepackage{epstopdf}
\usepackage{cleveref}

\startlocaldefs
 \newtheorem{theorem}{Theorem}
 \newtheorem{remark}{Remark}
  \newtheorem{proposition}{Proposition}
\endlocaldefs

\begin{document}

\begin{frontmatter}

\begin{fmbox}
\dochead{Research}


\title{Simulating brain rhythms using an ODE with stochastically
varying coefficients}


\author[
   addressref={aff1},                   
   email={benjamin.ambrosio@univ-lehavre.fr}   
]{\inits{B}\fnm{Benjamin} \snm{Ambrosio}}
\author[
   addressref={aff2, aff3, aff4},
   email={lsy@cims.nyu.edu}
]{\inits{L-S}\fnm{Lai-Sang} \snm{Young}}


\address[id=aff1]{
  \orgname{Normandie Univ, UNIHAVRE,LMAH,FR CNRS 3335,
ISCN}, 
  \postcode{76600},                                
  \city{Le Havre},                              
  \cny{France}                                    
}
\address[id=aff2]{%
  \orgname{Courant Institute of Mathematical Sciences, New York University},
  \city{New York},
   \postcode{NY 10012},
  \cny{USA}
}

\address[id=aff3]{%
  \orgname{Center for Neural Science, New York University},
  \city{New York},
  \postcode{NY 10003},
  \cny{USA}
}

\address[id=aff4]{%
  \orgname{Institute for Advanced Study},
  \city{Princeton},
  \postcode{NJ 08540},
  \cny{USA}
}


\begin{artnotes}
\end{artnotes}

\end{fmbox}


\begin{abstractbox}

\begin{abstract} 
The brain produces rhythms in a variety of frequency bands. Some are likely by-products of
neuronal processes; others are thought to be top-down. Produced entirely naturally, 
these rhythms have clearly recognizable beats, but they are very far from periodic
in the sense of mathematics. They produce signals that are broad-band, episodic,  
wandering in magnitude, in frequency and in phase; the rhythm comes and goes, degrading 
and regenerating. Rhythms with these characteristics do not match standard dynamical
systems paradigms of periodicity, quasi-periodicity, or periodic motion in the presence of
a Brownian noise. Thus far they have been satisfactorily reproduced only
 using networks of hundreds of integrate-and-fire neurons. In this paper, we tackle 
the mathematical question of whether signals with these properties
can be generated from simpler 
dynamical systems. Using an ODE with two variables inspired by the FitzHugh-Nagumo 
model, and varying randomly three parameters that control the magnitude, frequency and degree
of degradation, we were able to replicate the qualitative characteristics 
of these natural brain rhythms. Viewing the two variables as Excitatory and Inhibitory conductances
of a typical neuron in a local population, our model produces results that closely resemble 
gamma-band activity in real cortex, 
including the moment-to-moment balancing of E and I-currents seen
in experiments.

\parttitle{First part title} 
Text for this section.

\parttitle{Second part title} 
Text for this section.
\end{abstract}




\end{abstractbox}
%

\end{frontmatter}




\section{Introduction}  Rhythms, or oscillatory patterns of neural activity, occur ubiquitously
 in many parts of the central nervous system. One typically classifies 
them by their frequency bands. 
For example, $\alpha$-band activity  ($8$-$12$ Hz) is detected in the occipital lobes and is
associated with relaxed wakefulness, $\beta$-band rhythms ($12$-$30$ Hz) are related to muscles and movements, 
and  $\gamma$-rhythms ($30$-$90$ Hz) are implicated in  information 
transfer and are associated with cognitive processes. 
The origins and functional roles of these rhythms are active topics of current research, though much remains to be understood.

This paper is about simulating brain rhythms; the challenge is to generate mathematically signals that resemble the rhythms produced naturally by the brain. We will focus on gamma rhythms, and
for definiteness, we will base our study on gamma-band activity in the visual cortex, 
which has been the subject of detailed experimental studies e.g. \cite{Gra-1989,Hen-2005,Xin-2012}
and for which mechanisms are better understood (see below). See also the review article \cite{Car-2016}, 
and the modeling paper \cite{Cha-2018}.

Experimental data show that there are two aspects to the character of gamma rhythms: one is their oscillatory nature, 
and the other is their irregularity. In spite of their being called ``a rhythm", gamma rhythms are 
far from being periodic in the sense of mathematics. There is a recognizable beat, to be sure,
but spectral power density studies
show that gamma rhythms are broad-band, with wandering
frequencies and phases. Activity patterns are episodic; the beats are uneven in magnitude, 
degrading from time to time before the resumption of oscillatory behavior. 
One of the challenges to the dynamical systems community is to understand what kind of underlying dynamics can produce 
signals with these characteristics.  

Before explaining what is new, let us first review what has been done.
There have been a number of theoretical studies of gamma rhythms. An early model
is PING \cite{Whi-2000,Bor-2003}; similar models include \cite{Erm-1998,
Tie-2001} and others. These models produce highly regular population spikes, capturing
the oscillatory behavior of gamma rhythms well but not their irregular character. Other much cited works are \cite{Bru-1999, 
Bru-2000,Bru-2001}, which modeled gamma dynamics as periodic plus a noise term; that is not quite 
right still as it presumes a single intrinsic frequency. 
The irregular, episodic nature of gamma rhythms was captured in the multiple-firing events of 
\cite{Ran-2013} and further studied in \cite{Cha-2016,Cha-2018} using networks of hundreds to thousands of integrate-and-fire neurons. 

Our aim in this paper is not to build a biologically realistic model to reproduce the gamma activity in
real cortex. We seek to model only the phenomenon, or the signal, 
not the underlying biology, and to the extent that is possible, to do so analytically using as simple 
an equation as we can. 
The solution we came up with is defined by
a nonlinear ODE in two variables, with stochastically varying coefficients. 
In much the same way that the FitzHugh-Nagumo (FHN)
equations describe oscillatory behavior of a system with a voltage and a recovery variable, 
the two variables in our system can be thought of as describing the Excitatory and Inhibitory conductances produced
by a local neuronal population. We found that allowing the coefficients of our ODE to wander 
randomly produces results that best match the irregular character of gamma rhythms in real cortex.

The three main sections of this paper are Section 2, which 
describes the deterministic model, Section 3, where we introduce the stochasticity, and
Section 4, which examines the spectral properties and  frequencies of the signal produced,
to be compared to corresponding results in gamma-band activity in cortex.

\section{The deterministic model}

In Sect. 2.1 we present the model equations together with basic features of the dynamical
systems they generate. Sects. 2.2 and 2.3 study the dependence on the three parameters that are
allowed to vary.

\subsection{Equations and basic dynamical features}

 To motivate the equations we propose, we review quickly the {\it Recurrent Excitation-Inhibition} (REI) mechanism for gamma-band rhythms proposed in \cite{Cha-2018}. From the crossing of threshold by a few E-cells, recurrent 
excitation leads to the elevation of subthreshold activity and the spiking of more E and I-cells,
activating both populations essentially simultaneously. Depending on how the membrane potentials of 
neurons postsynaptic to the spiking cells are positioned, I-cells can step in to stop this developing 
event quickly, or the recruitment of more cells to join the event continues for a millisecond 
or two longer before the suppressive effects of I-spikes are felt. Since I-neurons are fairly densely connected to the rest of the local population, this pushback tends to hyperpolarize a good majority 
of the cells. The decay of the suppressive effect, and the depolarization of E-cells, 
helped along by the external drive, allows the scenario to be repeated again. This explains 
the production of a rhythm. Biophysical time constants cause the cycles to have 
frequencies in the gamma band.

In our model $u$ and $v$ represent the absolute values of the E and I-conductances
of a typical neuron in  a local population.
We propose that their dynamics are described by the following equations: 
\begin{equation}
\label{eq:a-y}
\left\{\begin{array}{rcl}
\epsilon u_t&=&u(-K(u-a_1)(u-a_2)-v)\\ 
v_t&=&\gamma v(bu-v+c)
\end{array}
\right.
\end{equation}
where  $a_1, a_2, b$ and $c$ are fixed parameters with 
\begin{equation}
\label{eq:parameters-ini}
\begin{array}{c}
 a_1=-0.01, \ \ a_2=0.1, \ \ b= 11.9, \ \ c=6.6 \times 10^{-4},
\end{array}
\end{equation}
while $\epsilon, \gamma$ and $K$ are parameters taking values, initially at least, in 
$$\epsilon \in [0.01,1], \ \ \gamma \in [1,25] \quad \mbox{and} \quad K\in [30,100]\ .
$$ 
These equations are inspired by the FHN system, which also has a slow-fast structure.
The cubic nullcline in FHN is replaced by a parabola in the first equation. Factors of $u$ and $v$ are added respectively in front of the first and second equations to constrain the solutions to the upper-right quadrant as well as to better fit conductance oscillations. These modifications are also inspired by the Leslie-Gower model, see \cite{Amb-2018}. 

In the rest of this paper, we will adopt the following notation:
\[F(u,v)=u(-K(u-a_1)(u-a_2)-v),\,\ \ G(u,v)=\gamma v(bu-v+c)\]
\[f(u)=-K(u-a_1)(u-a_2),\,\ \ g(u)=(bu+c)\ .\]

The nullclines of Eqns (\ref{eq:a-y}) are then given by
\[u=0, \, v=f(u)\]
for the first equation, and
\[v=0, v=g(u)\]
for the second equation. Note that the  polynomial $f(u)=-K(u-a_1)(u-a_2)$ reaches its maximum for $u=0.5(a_1+a_2)$ with a value of $0.25K(a_2-a_1)^2$.

Our choice of parameters allows only one intersection between the quadratic function $f(u)$ 
and the linear function $g(u)$ in the positive quadrant. This is ensured by the condition:
\[b\frac{a_1+a_2}{2}+c>0.25(a_2-a_1)^2K\]
which gives
\[K<176.8.\]
There are four stationary points in the positive quadrant, the region of interest. They are
\[(0,0), (0,c), (a_2,0), (u^*,v^*)\]
where $u^*$ is the positive solution of
\[K(u-a_1)(u-a_2)+bu+c=0 \]
and
\[v^*=bu^*+c.\]
See Figure 1, which gives a sense of the global dynamics and basic structures
of the system.

The next result confirms that the region depicted in Figure 1 is indeed positively invariant.

\begin{theorem}
\label{th:lc}
The positive $(u,v)$-quadrant is invariant under the dynamics defined by Eq. (\ref{eq:a-y}), 
and there exists a bounded absorbing set to which all solutions enter.
\end{theorem}

\begin{figure}
\begin{center}
 \label{fig:FP}
\includegraphics[height=9cm, width=12cm]{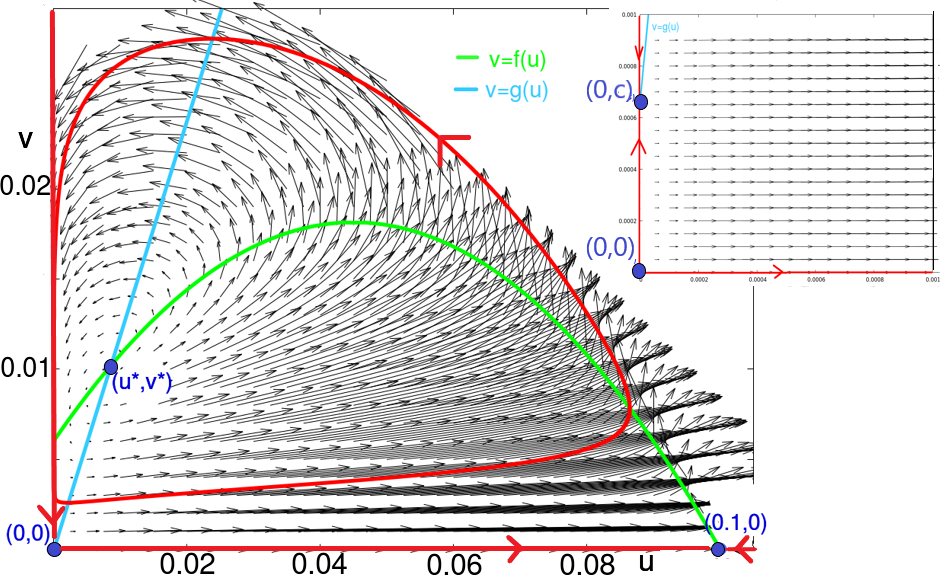}

\end{center}
\caption{This figure shows the fixed points, nullclines, and vector field, as well as the
trajectories lying in the sets $u=0$, $v=0$ and the limit-cycle for $K=60$ and $\epsilon=0.1$. For these values of parameter $ (a_2,0)$ is a saddle and $(u^*,v^*)$ is a source.  The two fixed points 
$(0,0)$ and $(0,c)$ are too close to be discernible in this figure. In the inset, 
we zoom in to visualize them: $(0,0)$  is a source, and $(0,c)$  is a saddle.
}
\end{figure}

\begin{proof}
The fact that the positive quadrant is positively invariant follows from the fact that $u=0, v_t=\gamma v(-v+c)$ and $\epsilon u_t=-u(K(u-a_1)(u-a_2)), v=0$, are solutions lying on the $u$ and $v$-axes. For the existence of an absorbing set, we compute 
\[\frac{d}{dt}(\epsilon u^2+v^2).\]
We find:
\begin{equation}
\label{eq:absset}
\frac{d}{dt}(\epsilon u^2+v^2)=2\big[-Ku^4+K(a_1+a_2)u^3-Ka_1a_2u^2-u^2v+b\gamma uv^2-\gamma v^3+c\gamma v^2   \big]
\end{equation}
By using Young inequality, for $u,v >0$,
\[uv\leq \frac{u^p}{p}+\frac{v^q}{q}, \, \frac{1}{p}+\frac{1}{q}=1 .\]
We can split the right hand side of \eqref{eq:absset} into a  polynomial in $u$ of degree $4$, and a polynomial in $v$ of degree $3$, with both negative leading coefficients. Since the solution lies in the positive quadrant, by polynomial comparison, we obtain:
\begin{equation}
\label{eq:ieqas}
\frac{d}{dt}(\epsilon u^2+v^2)\leq -A(\epsilon u^2+v^2)+B
\end{equation}
where $A>0$ and $B>0$ can be chosen independently of initial conditions. Integrating \eqref{eq:ieqas} leads to the existence of an absorbing set attracting all trajectories.
\end{proof}

\subsection{Varying the parameter $\epsilon$}

\begin{figure}
 \label{fig:eps=0001-10K=60}
\includegraphics[height=5cm, width=6cm]{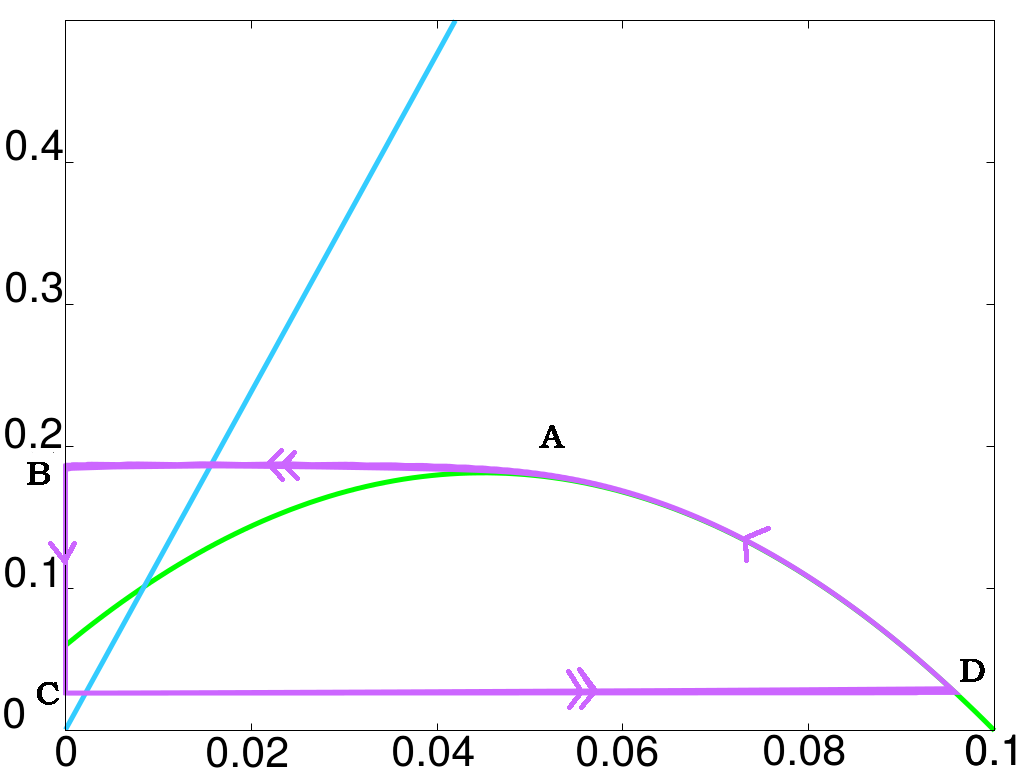}
\includegraphics[height=5cm, width=6cm]{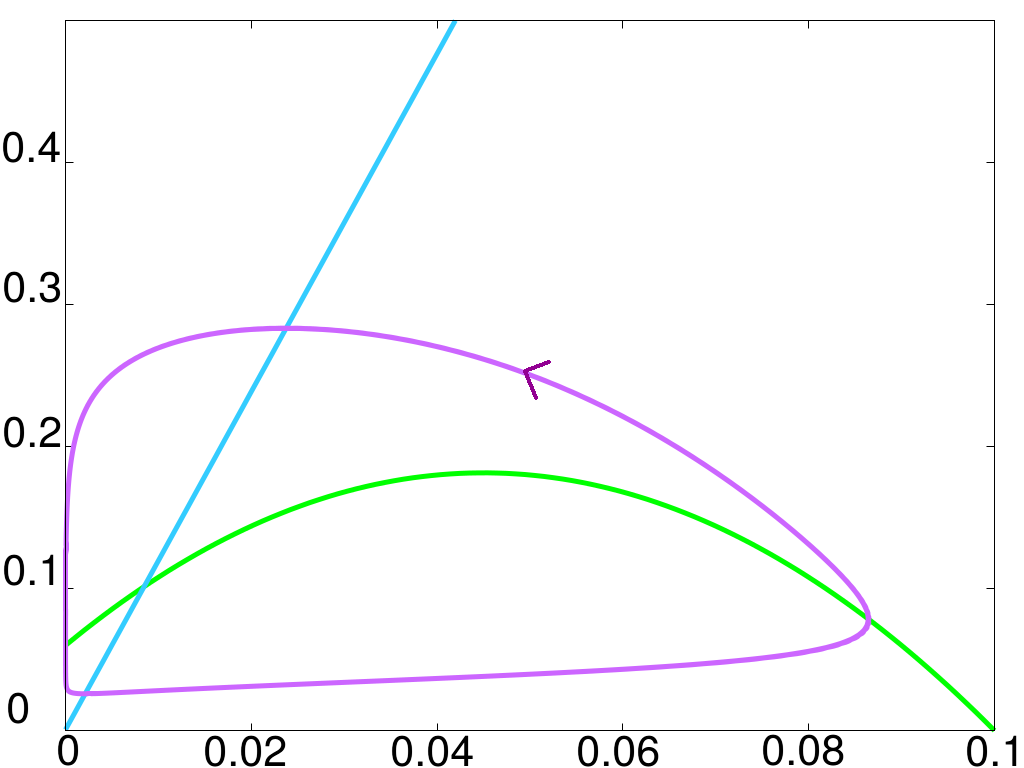}\\
\includegraphics[height=5cm, width=6cm]{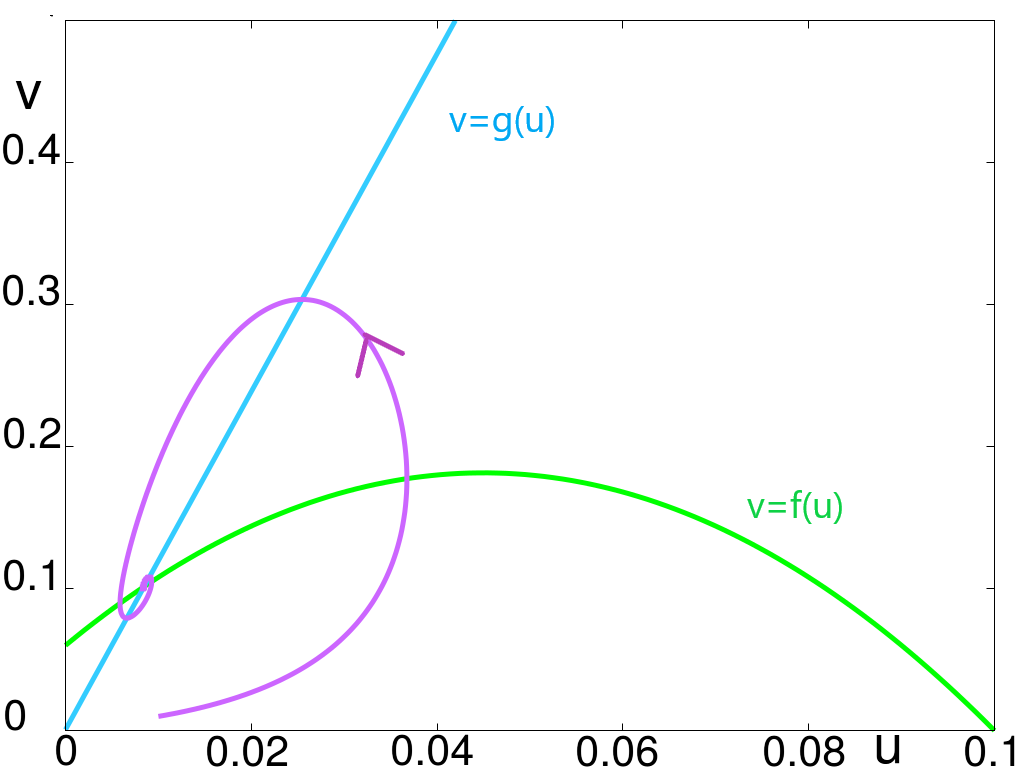}
\includegraphics[height=5cm, width=6cm]{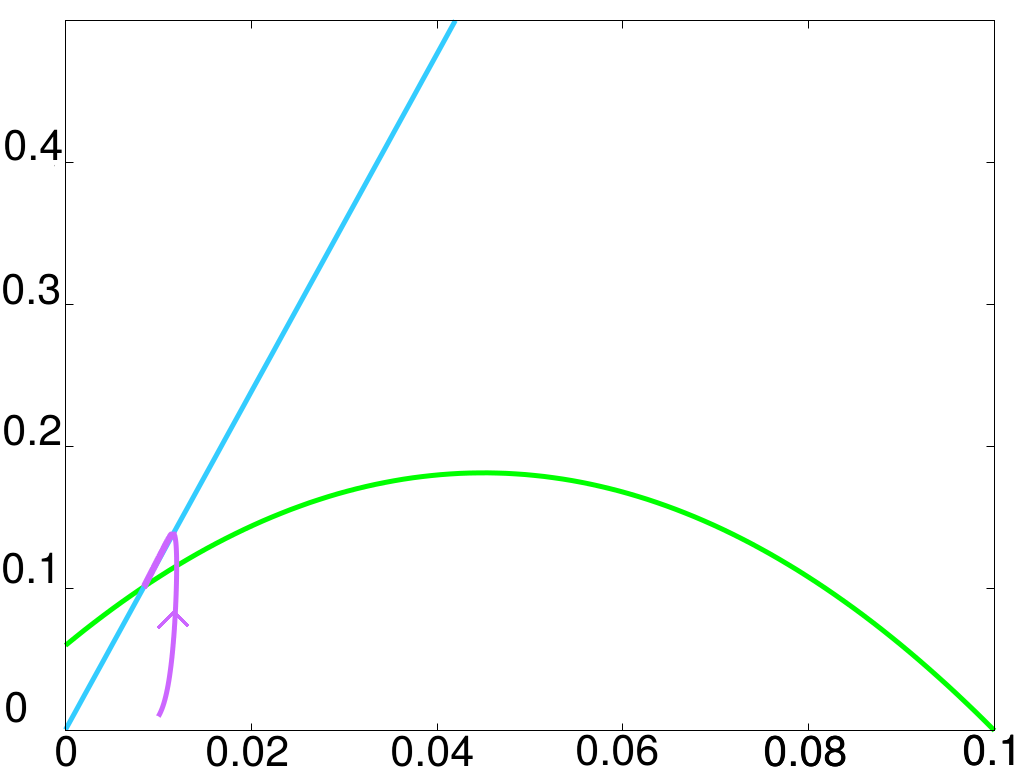}
\label{fig:eps=1K=60}
\caption{This figure illustrates the effect of varying $\epsilon$, from small to large values. The parameter $K$ is set to $K=60$. The top left panel illustrates $\epsilon=0.001$, the top right, $\epsilon=0.1$, bottom left, $\epsilon=1$ and bottom right $\epsilon=10$.}
\end{figure}

Recall that we have three parameters: $\epsilon, \gamma$ and $K$. We first explain the role
of $\epsilon$, fixing for now $\gamma=1$, and studying the dynamics as $\epsilon$ is varied
for each value of $K$. The dynamics of the system from $\epsilon$ very small to very large
for $K=60$ (a fairly typical value of $K$) are summarized in Figure 2. 
For $\epsilon \ll 1$, Eq. (\ref{eq:a-y}) describes a slow-fast system
with a limit cycle as can be seen in the top two panels of Fig 2. As $\epsilon$ is increased,
this limit cycle turns into a sink somewhere between $\epsilon=0.1$ and $\epsilon=1$. 
For large $\epsilon$, e.g at $\epsilon=10$, one can show that the the critical attractive 
manifold is the line $\Delta: v=g(u)$. All the trajectories reach this line fast and then follow it slowly  
toward the fixed point $(u^*,v^*)$ as can be seen in the bottom right panel.

Below we will give the analysis for $\epsilon$ very small, as well as the Hopf bifurcation that takes
the limit cycle to the sink.

\subsubsection{The case of $\epsilon \ll 1$}

For $\epsilon$ small enough a slow-fast analysis allows to compute the limit-cycle up to an $O(\epsilon)$ order. The behavior can be described geometrically as follows. We denote by $\mathcal{C}$ the curve  $v=f(u)$. For $\epsilon$ small enough, a trajectory starting from the right side of $\mathcal{C}$ will increase along the curve ($v_t>0$ there) until it reaches the maximum point $A=(\frac{a_1+a_2}{2},f(\frac{a_1+a_2}{2}))$. This is a jump point, see \cite{Kru-2001}, i.e., from there the trajectory leaves $\mathcal{C}$ and goes 
at high speed to reach a neighborhood of the point $B=(0,f(\frac{a_1+a_2}{2}))$. After that, since at first $u_t<0$, the trajectory remains stuck near the line $u=0$. It goes down ($v_t<0$) until it crosses the point $(0,f(0))=(0,-Ka_1a_2)$, at which point $u_t$ becomes positive. This is a fold point but not a jump point, see \cite{Kru-2001}.  Dynamics near this point have been analyzed in \cite{Amb-2018}. See also \cite{Wan-2019} and references therein cited. The trajectory continues to follow the axis $u=0$ until it reaches a point $C$ on the axis $u=0$ which is significantly below $(0,f(0))$.  Here there is the possibility of the so-called \textit{canard} phenomenon,  see \cite{Ben-1981,Kru-2001,Szm-2001}.
At $C$, the trajectory leaves the axis $u=0$ and goes very quickly toward the point $D$ on $\mathcal{C}$ with the same ordinate as $C$. This gives a qualitative description of the limit-cycle. For $\epsilon$
sufficiently small, precise statements can be rigorously deduced from Geometrical Singular Perturbation  Theory. Good reviews can be found in \cite{Hek-2010,Jon-1995,Kap-1999,Kru-2001}. 

Let $\Gamma'$ be the closed curve defined by:
\[\Gamma'=[A,B]\cup [B,C]\cup [C,D] \cup \zeta \]
where, $\zeta \subset \mathcal{C}$ is the arc from $D$ to $A$.   
\begin{theorem}
\label{th:eps-small}
For $\epsilon>0$ sufficiently small, there is a limit cycle $\Gamma$ within distance
$O(\epsilon)$ of $\Gamma'$.
\end{theorem}  

\begin{remark}
We did not consider the uniqueness of the limit-cycle here. For a proof in the case $\epsilon$ small, see \cite{Wan-2019}.
\end{remark}

We point out that the system defined by \eqref{eq:a-y} provides a simple example, in a Neuroscience context, in which canard solutions emerge and can be computed explicitly. Thanks to the polynomial expression of  the vector field, the computations performed in \cite{Amb-2018}  become simpler and explicit around the point $(0,f(0))$. Details are given in the Appendix for the convenience of the reader.

\begin{figure}
 \label{fig:FromSinkToLC}
 \includegraphics[height=4cm, width=5cm]{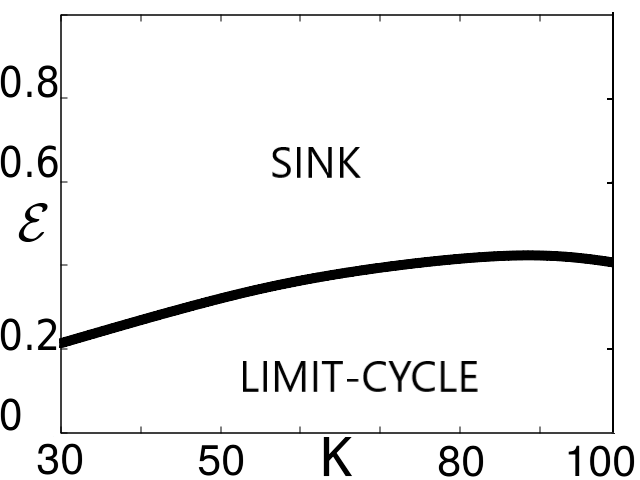}\\
\includegraphics[height=4cm, width=4cm]{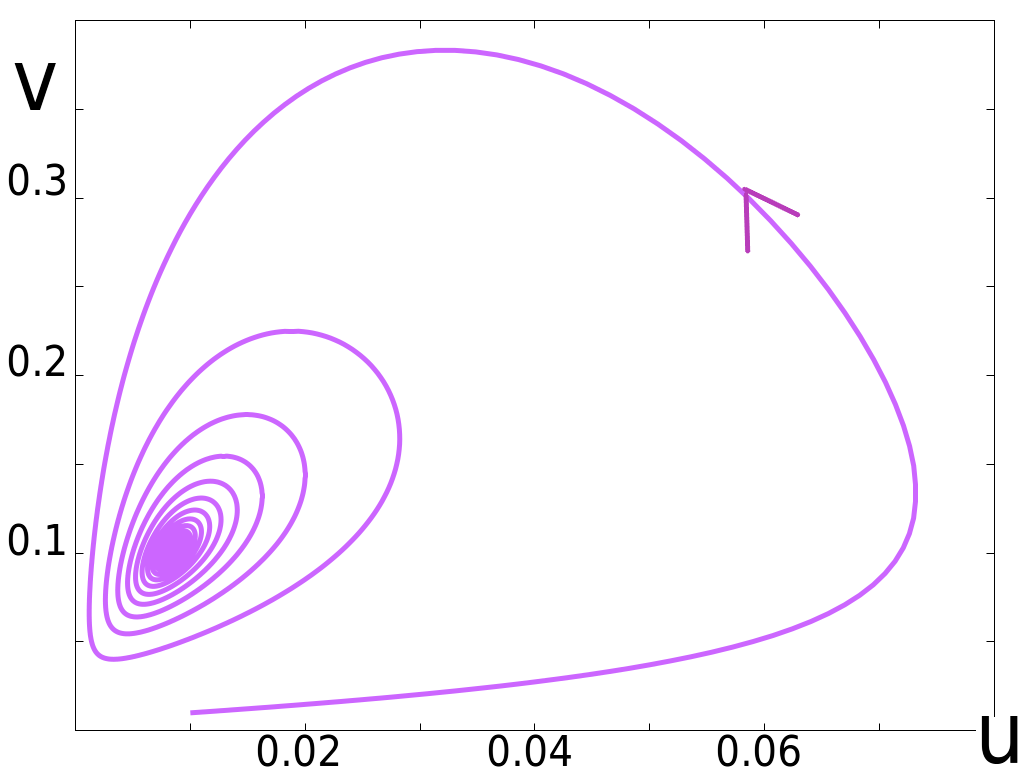}
\includegraphics[height=4cm, width=4cm]{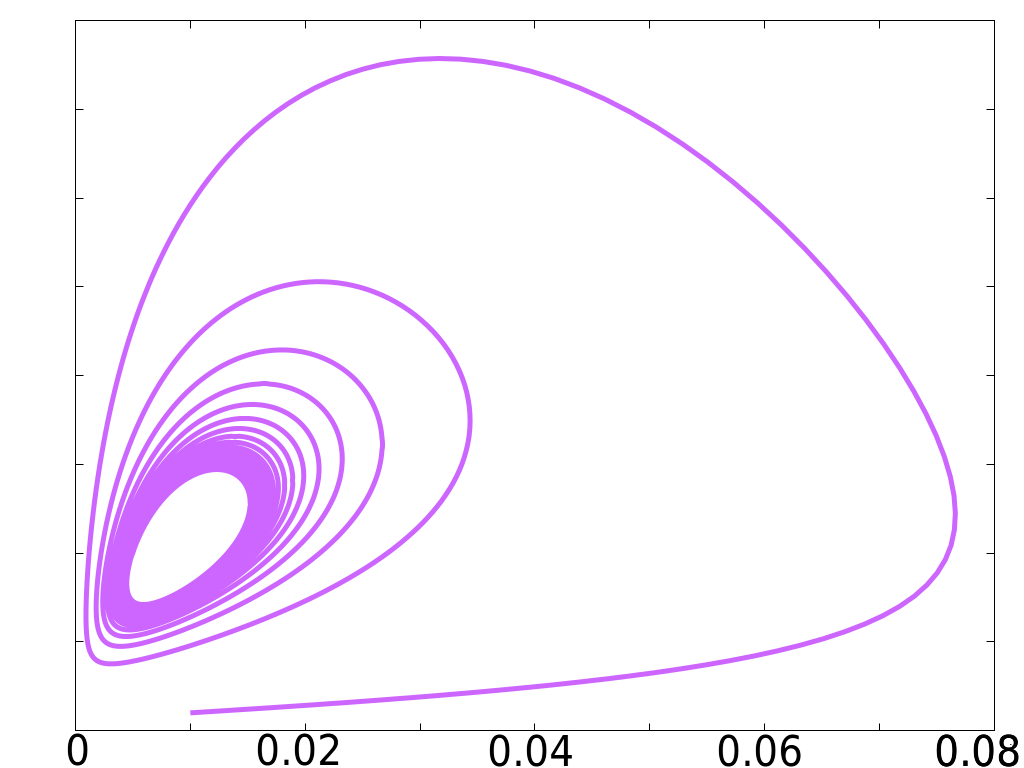}
\includegraphics[height=4cm, width=4cm]{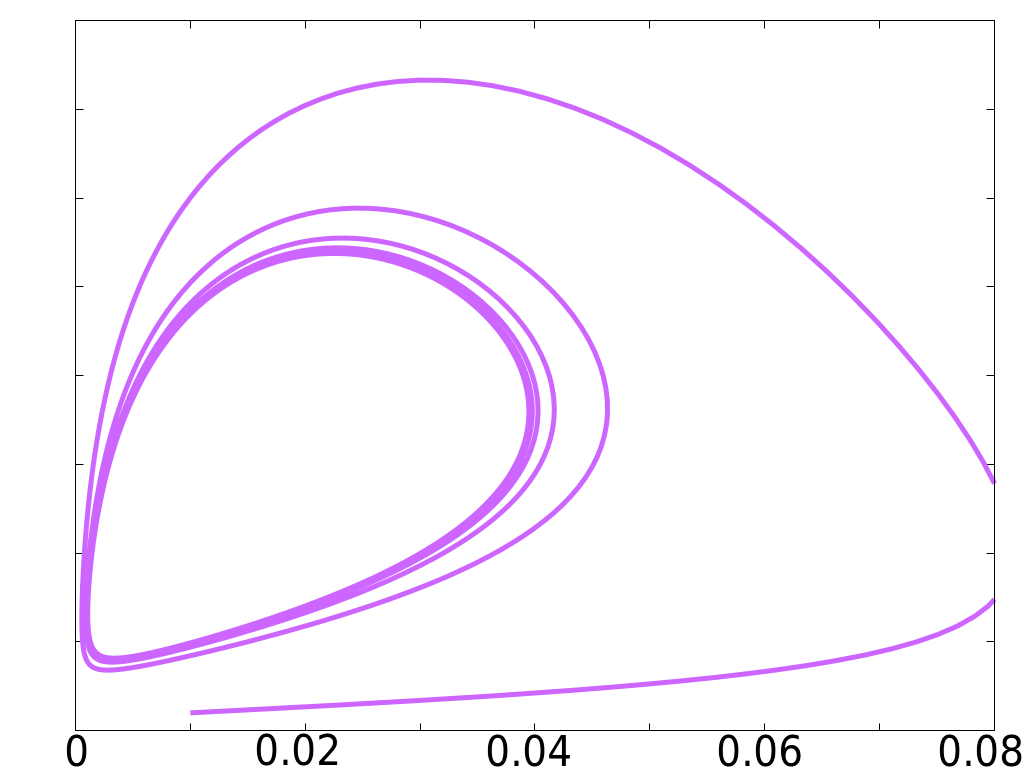}
\caption{Illustration of the Hopf bifurcation. In the top panel we have plotted the Hopf bifurcation diagram in the $K,\epsilon$ plane. The bottom panels illustrate the bifurcation for $K=60$ as $\epsilon$ decreases. 
Left: $\epsilon=0.4$, a trajectory spiraling towards a sink. Middle: $\epsilon=0.36$, trajectories
accumulating on a limit cycle following the sink's loss of stability. Right: $\epsilon=0.3$, the limit cycle
growing in size.}
\end{figure}

\subsubsection{Hopf bifurcations}

As indicated in Sect. 2.1, the range of $\epsilon$ of interest  is $[0.1, 1]$, and it is in this range
of $\epsilon$ that the limit cycle turns into a sink as shown in Fig 2. We now give more detail 
on this bifurcation, specifically the Hopf bifurcation that occurs at $(u^*, v^*)$ where $(u^*, v^*)$ 
is the unique fixed point in the interior of the positive quadrant; see Sect. 2.1.

The Jacobian matrix at fixed points is
\[J=\begin{pmatrix}
\frac{1}{\epsilon}(-3Ku^2+2K(a_1+a_2)u-Ka_1a_2-v)& -\frac{1}{\epsilon}u\\
bv & -2v+(bu+c)\ 
\end{pmatrix}.\]
Substituting in $v^*=bu^*+c$, we obtain at $(u^*,v^*)$, that
\[J_{(u^*,v^*)}=J^*=\begin{pmatrix}
-\frac{1}{\epsilon}Ku^*(2u^*-(a_1+a_2))& -\frac{1}{\epsilon}u^*\\
b(bu^*+c) & -(bu^*+c)
\end{pmatrix},\]
which gives
\[\det(J^*)=\frac{1}{\epsilon}u^*(bu^*+c)(K(2u^*-(a_1+a_2))+b)\]
while 
\[tr(J^*) = -\frac{1}{\epsilon}Ku^*(2u^*-(a_1+a_2)) - (bu^*+c)\ .\]
From the above expressions, we deduce the following proposition.

\begin{proposition}
For $K\in [30,100], det(J^*)>0$. It follows that for each $K$ there exists a value of $\epsilon$ at which
a Hopf bifurcation occurs. This value is given by:
\[\epsilon=Ku^*\frac{a_1+a_2-2u^*}{bu^*+c}.\] 
\end{proposition}

We close this section by an application of the Poincare-Bendixon theorem to our system.
\begin{theorem}
Each trajectory starting in the region $\{u>0,v>0\}$ either converges to $(u^*,v^*)$ or evolves towards a limit-cycle. For $\epsilon<Ku^*\frac{a_1+a_2-2u^*}{bu^*+c}$, it converges towards a limit-cycle.
\end{theorem}
\begin{proof}
The proof follows from the analysis of the nullclines and the nature of fixed points. 
\end{proof}

\subsection{Dependence of dynamics on the parameters $\epsilon, K$ and $\gamma$}

Continuing to keep $\gamma=1$, we first examine the dynamics of 
Eq. (\ref{eq:a-y}) as functions of $K$ and $\epsilon$. Simulation results are shown in Fig 4.
Notice first that these results are consistent with those in Fig 3 with regard to increasing
$\epsilon$ for fixed $K$. What is new here is the effect of varying $K$ for each $\epsilon$. 
Fig 4 shows clearly that larger $K$ corresponds to larger excursions by $u$ and $v$. 
This means 
\begin{itemize}
\item[(i)] when solutions are attracted to a limit cycle, the limit cycle has a larger diameter for larger $K$; and 
\item[(ii)] whether the $\omega$-limit set is a  limit cycle or a sink, it is located closer to $u=0, v=0$ 
for smaller values of $K$.
\end{itemize}

\begin{figure}
 \label{fig:eps=01-04K=30-100}
\includegraphics[height=4cm, width=5cm]{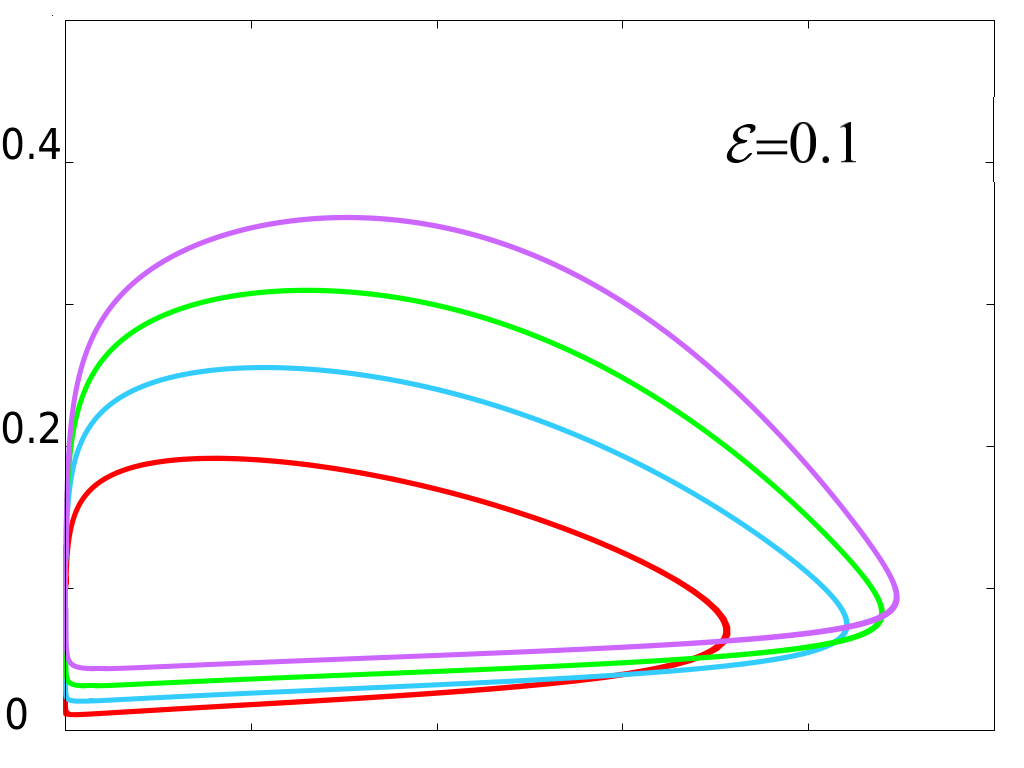}
\includegraphics[height=4cm, width=5cm]{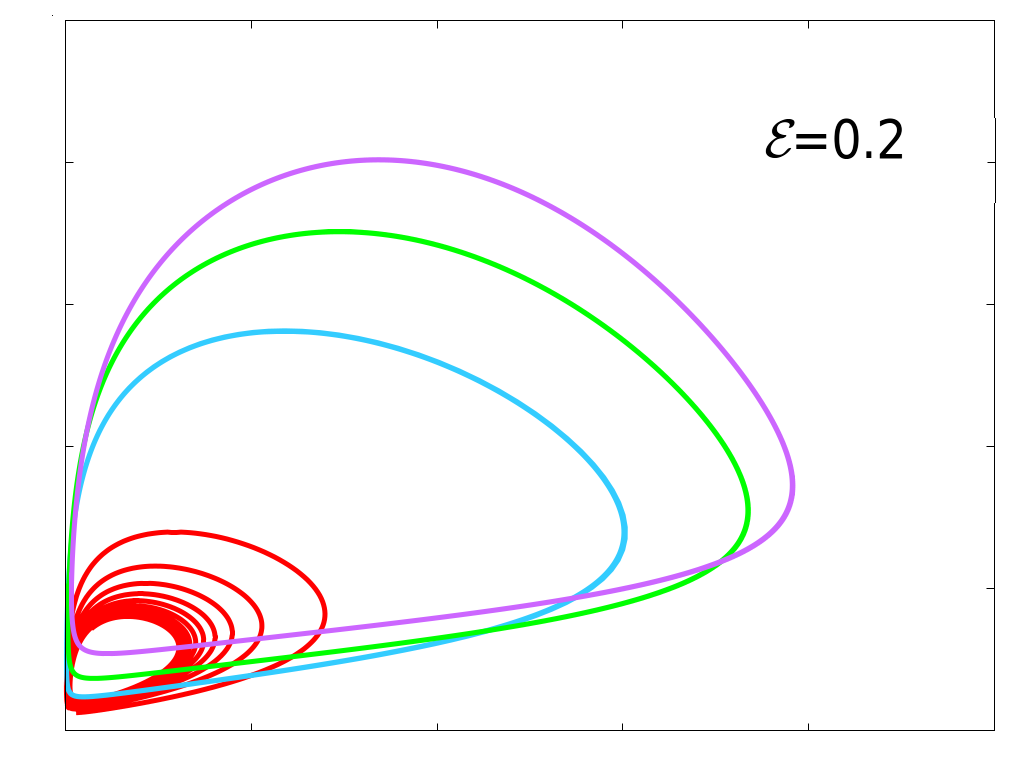}\\
\includegraphics[height=4cm, width=5cm]{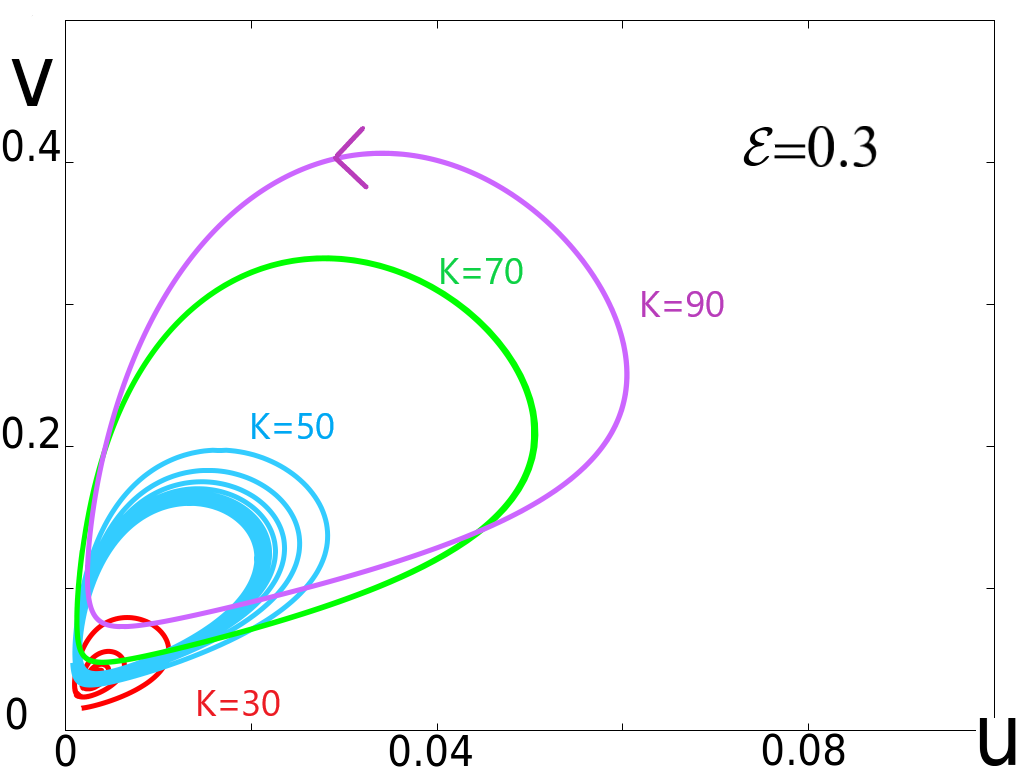}
\includegraphics[height=4cm, width=5cm]{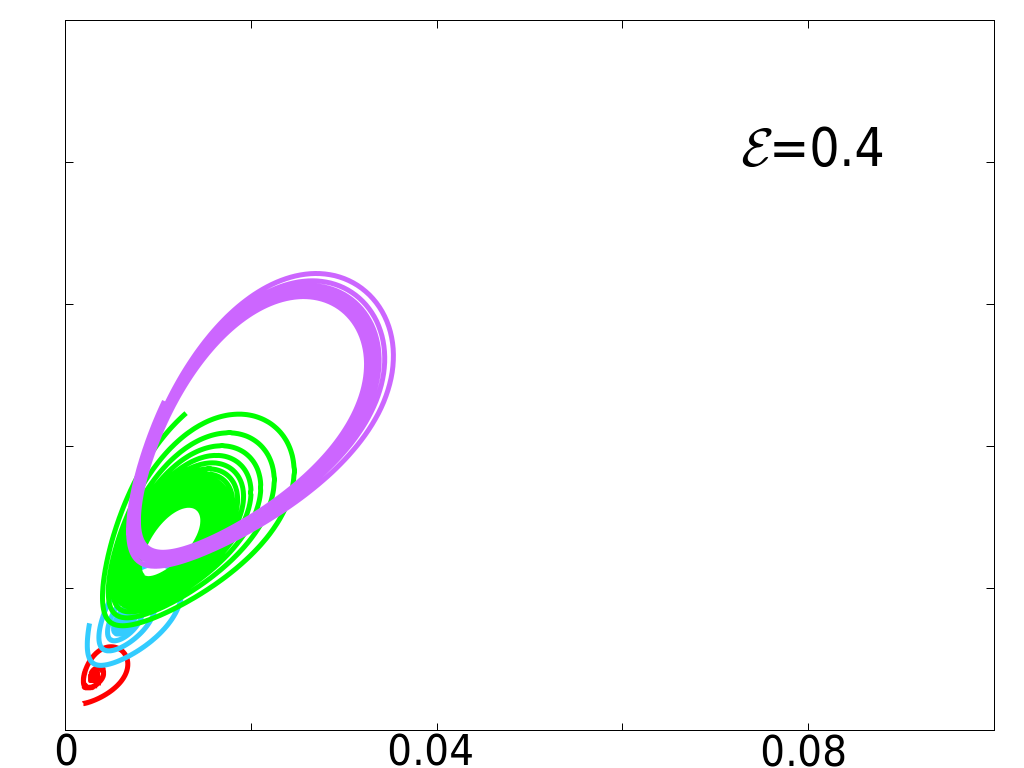}
\caption{This figure gives a panorama of sample trajectories within the paremeters range of interest. 
Four panels corresponding to $\epsilon=0.1, 0.2, 0.3$ and $0.4$ are shown. In each panel,
trajectories for different values of $K$ are depicted by different colors: 
$K=30$ (red), $50$ (cyan), $70$ (green) and $90$ (purple).}
\end{figure}

\medskip
Finally, we examine the effect of varying $\gamma$. From the equations, it is clear that trajectories
of Eq. (\ref{eq:a-y}) will trace out the same curves as long as $\gamma\epsilon$ remains constant;
and that varying $\gamma$ keeping $\gamma\epsilon$ fixed corresponds to changing the speed with which one moves along these curves. For example, at $K=60$, for values of $\epsilon=0.1$
and $\gamma=1$,  numerical simulation gives a limit cycle with period $\sim 44$  ms (equivalently
frequency around 22 Hz). For $\epsilon=0.01, \gamma=10$, the period becomes $4.4$ ms (frequency around 225 Hz). 

\begin{proposition}
\label{prop:freq} For each fixed $K$, the curves traced out by the trajectories of Eq. (\ref{eq:a-y}) 
depend only on  $\epsilon \gamma$. Fixing $K$ and $\epsilon \gamma$, and varying $\gamma$,
velocities are proportional to $\gamma$; in particular, the frequency of the limit cycle is proportional to $\gamma^{-1}$.
\end{proposition}

The meaning and main general effects of the variation of parameters $\epsilon, K$ and $\gamma$ in Eq. (\ref{eq:a-y})
can be summarized as follows:
{\it \begin{itemize}
\item increasing $\epsilon$ changes the dynamical regime from one with a limit cycle 
in a slow-fast 
system to one with an attractive fixed point; 
\item $K$ controls the sizes of the excursion of $(u,v)$ in the system's oscillatory behavior: 
in general, the larger $K$, the larger the excursions; while
\item 
 for each fixed value of $\epsilon \gamma$, the magnitude of $\gamma$ controls the frequencies 
of the limit cycle.
\end{itemize}}

\medskip
As we will show momentarily, these are the parameters we need to vary to produce the irregularities
seen in gamma rhythms.

\section{The stochastic model}

As discussed in the Introduction, there are two facets to gamma rhythms as observed
in real cortex: one is their oscillatory nature; the other is their irregular, episodic character.
The deterministic system in Section 2 provided the underlying oscillations. Here we create
irregularity by adding randomness to the deterministic model. 
We begin by explaining the motivation for what is done.

A standard way to model random dynamical systems is to use a stochastic differential equation (SDE).
These ideas were used in \cite{Bru-1999,Bru-2000}, who studied networks of sparsely coupled
integrate-and-fire neurons, focusing on the limit where system size tends to infinity
and the number of connections is infinitely small compared to system size. Arguing that distinct
neurons are likely to have disjoint sets of presynaptic cells, 
the authors of 
\cite{Bru-1999,Bru-2000} modeled neuronal dynamics by an equation consisting of
a deterministic part describing meanfield activity plus a Gaussian noise that is independent 
from neuron to neuron, and gamma rhythms were modeled as regimes following a supercritical 
Hopf bifurcation. Another relevant paper is \cite{Bru-2003}. Here the authors assumed that gamma rhythms consisted of random noise 
on top of a periodic motion; they were primarily interested in the periodicity. 

Since the publication of these theoretical results, new experimental data were obtained, 
 including much higher
connection probabilities between neurons \cite{Hol-2003,Osw-2011} and 
the broad-band nature of gamma rhythms in real cortex \cite{Hen-2005}. 
The latter in particular suggests that gamma activity will likely fit better with 
dynamics that are not tied to a single period, something assumed in previous works. 
As to the nature of the irregularity, it was shown in real data that 
gamma power and frequencies {\it wander} \cite{Xin-2012}, i.e. these quantities drift, 
but similar gamma patterns are often repeated for tens, sometimes up to one or two hundred, 
milliseconds. 
The same was observed in simulations 
using networks of integrate-and-fire neurons \cite{Cha-2018}, and 
the phenomenon can be explained by the 
REI mechanism (see Sect. 2.1): In the aftermath of a large firing event that involves a larger than usual 
fraction of the E and I-population, inhibitory conductances  tend to 
be strong. Disinhibition occurs at roughly the same time for many neurons, 
increasing the likelihood that another major event will be precipitated when the march towards 
threshold is too synchronized. Similarly, when the rhythm degrades, such as when the nucleation 
of a firing event is interrupted too quickly by I-firing, the situation will likely remain ambiguous
for some time before another large event can occur.

To model the wandering character of gamma power and frequency, as well as the occasional
degradation of the rhythm, we propose to use an ODE whose coefficients are not fixed but 
are allowed to drift freely and
randomly -- effectively performing random walks -- within certain designated ranges.

 In more detail, we first specify parameter ranges $[K_{\min}, K_{\max}]$,
$[\epsilon_{\min},\epsilon_{\max}]$ and $[f_{\min},f_{\max}]$ for $K, \epsilon$ and $\epsilon \gamma$
respectively. (For the simulations shown in Figs 5 and 6, we used $[K_{\min}, K_{\max}]=[30,100]$, 
$[\epsilon_{\min},\epsilon_{\max}]=[0.04,0.1]$
and $[f_{\min},f_{\max}]=[0.2, 0.5]$.) Let $\mathcal U^1_i, \mathcal U^2_i, \mathcal U^3_i, i=1,2, \cdots,$ be 
independent random variables uniformly 
distributed on $[-1,1]$. 
Starting from initial values of $K, \epsilon$ and $\gamma$ within
the specified ranges, we update these parameters every $0.1$ ms. 
At the $i$th step, we let 
\[K=K(1+0.1\mathcal{U}^1_i), \]
constraining $K$ to $[K_{\min}, K_{\max}]$ according to the rule that if
$\mathcal U^1_i=u$ and $K(1+0.1u)$ falls outside of
$[K_{\min}, K_{\max}]$, then we set $K=K(1-0.1u)$. Next, we update $\epsilon$ 
by letting 
\[\epsilon=\epsilon+0.01\mathcal{U}^2_i\ \]
constraining $\epsilon$ to $[\epsilon_{\min},\epsilon_{\max}]$ as before. 
Finally, we set
\[\gamma=\gamma+0.1\mathcal{U}^3_i\]
if $\epsilon \gamma \in [f_{\min},f_{\max}]$. If not, we redraw $\mathcal U^3_i$ until $\epsilon \gamma \in [f_{\min},f_{\max}]$.

 Recall that it is the product $\epsilon \gamma$ that determines the curves
traced out by the trajectories of the system  (Proposition 2), 
and $\epsilon \gamma \in [0.2, 0.5]$
corresponds to $\epsilon \in [0.2, 0.5]$ in Figs 2--4, where $\gamma$ was set $=1$.  
Thus to simulate gamma rhythms, the parameters above 
are chosen so that most but not all of the time, 
the dynamics are oscillatory. Once parameters that
produce suitable qualitative behaviors are located, it is generally simpler to adjust the values
of $u, v$  or the mean frequencies of the oscillations by modifying slightly the two equations of Eq. \eqref{eq:a-y} (e.g. by inserting a scaling coefficient in front).

\begin{figure}[!htb]
 \label{fig:uvsto}
\includegraphics[height=3.5cm, width=12cm]{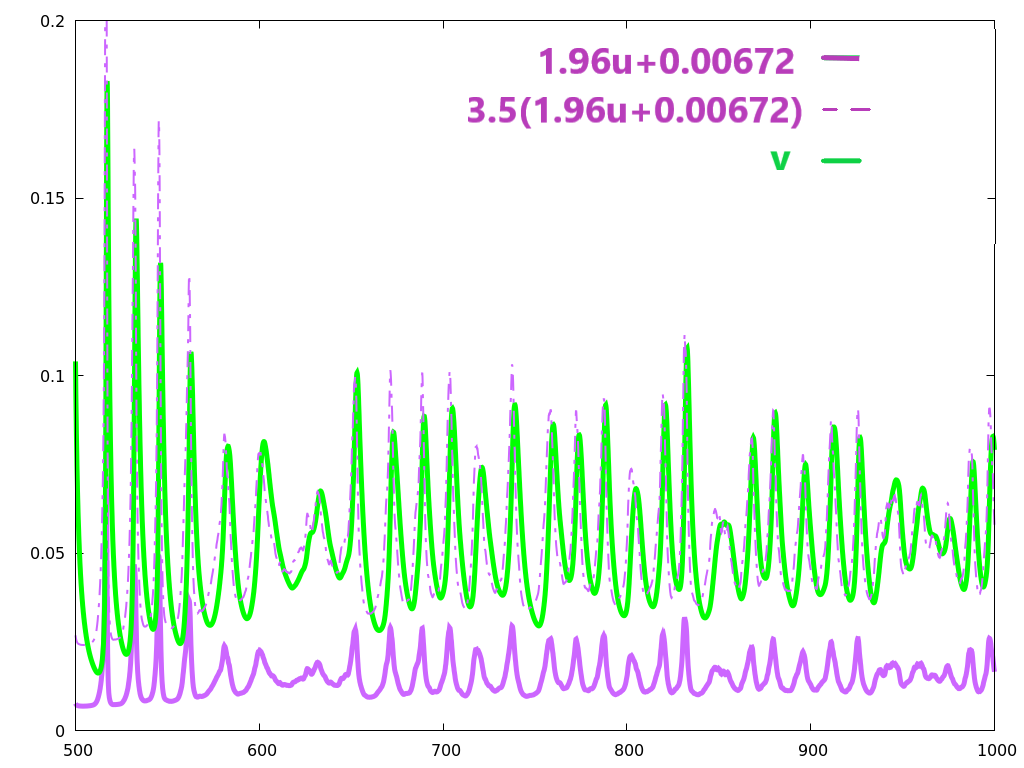}\\
\includegraphics[height=3.5cm, width=12cm]{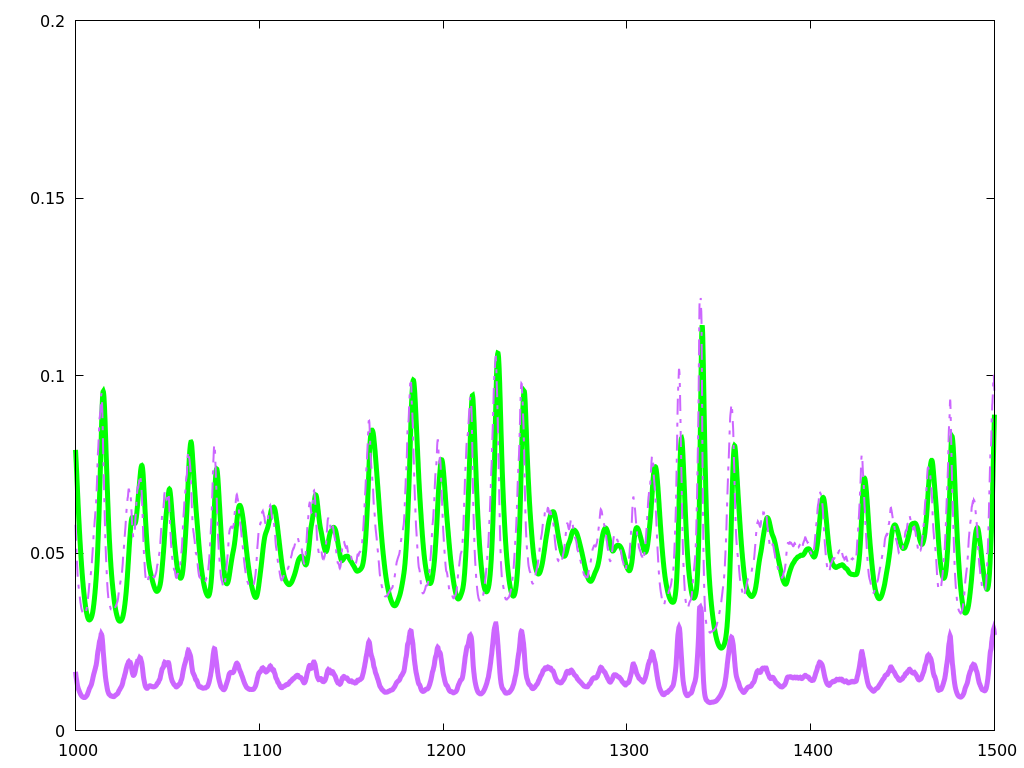}\\
\includegraphics[height=3.5cm, width=12cm]{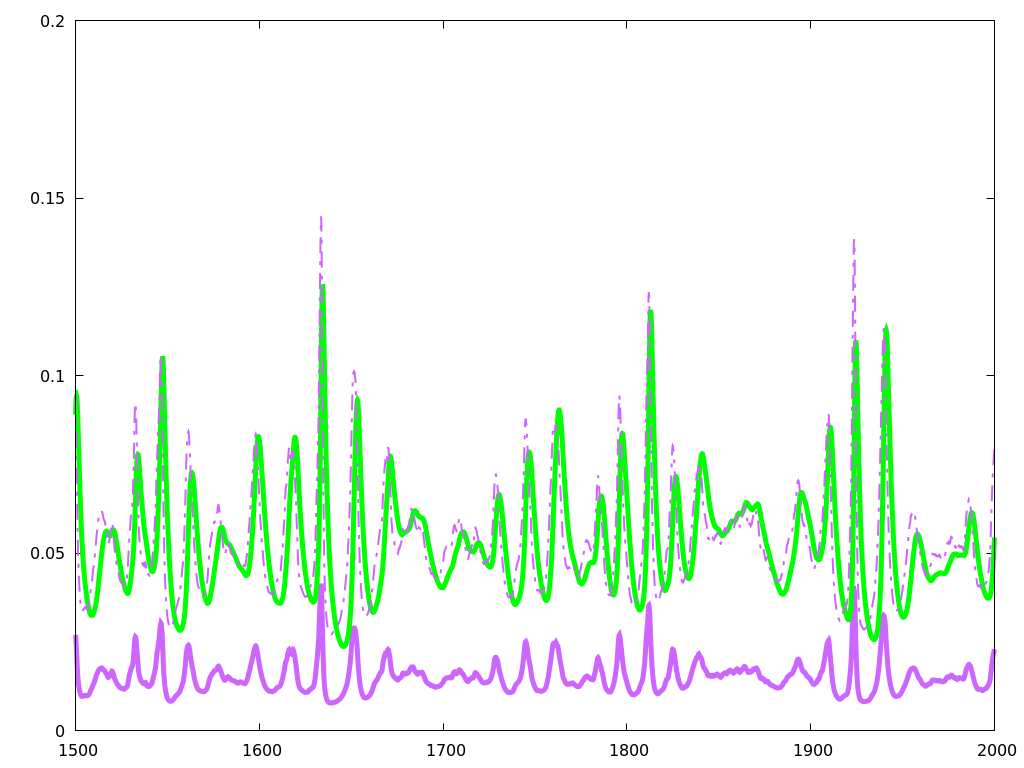}\\
\includegraphics[height=3.5cm, width=12cm]{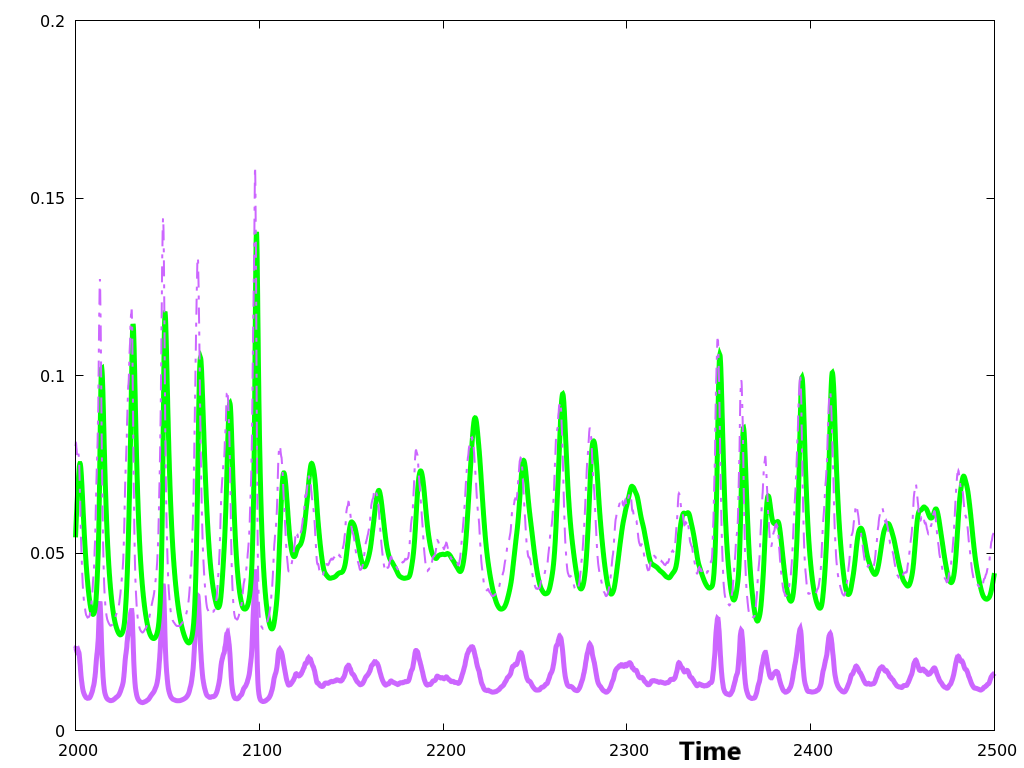}\\
\caption{This figure represents the evolution of the stochastic version of the system (\ref{eq:a-y}). 
The time windows in the four panels are (in ms) [500-1000],[1000-1500],[1500-2000] and [2000-2500]. 
The parameters used are $K \in [30,50], $ $\epsilon \in [0.04, 0.1]$, $\gamma \epsilon \in  [0.2,0.5]$. 
We plotted 
$\bar u =1.96u+0.00672$ in solid purple, to be thought of as representing E-conductance, 
$v$ in green, representing I-conductance. Since the ratio of E-current to E-conductance
is roughly three to four times that of the ratio of I-current to I-conductance, we have plotted
also $3.5 \bar u$ in dashed purple. Note the tight relationship between $3.5 \bar u$ and $v$.}
\end{figure} 

 Figure 5 shows a solution of the stochastic version of Eq. \eqref{eq:a-y} 
adjusted as indicated to simulate the E and I-conductances
of a typical neuron in a local population under drive thus producing gamma rhythms.
The irregular, episodic nature of the rhythm is 
clearly visible. By definition, the E-current entering a neuron is defined to be its E-conductance
times a factor proportional to the distance of membrane potential 
to the E-reversal potential, and the same is true for I-currents. As this factor for E-current is 
$3$ to $4$ times that for I, 
we have also plotted (in dash) a graph that is 
$3.5$ times the height of the E-conductance. Modulo a multiplicative constant, then, 
the dashed purple and green plots can be thought of as approximations of
 E and I-currents respectively.
 
We remark on the tightness with which the green plots (I-current) follow the dashed purple plots
(E-current). There is a well known theory of balanced states \cite{Vre-1998} 
that asserts that in the limit as system size tends to infinity, E-currents and I-currents are balanced
when averaged over time. But the experimental results of \cite{Oku-2008} and subsequently
the modeling paper \cite{Cha-2018} show that much more than that is true, namely that these
currents are in fact roughly balanced from moment to moment, not just when averaged over time.
This is entirely consistent with the REI mechanism, which tells us that in most
firing events, both E and I-neurons are recruited;  the size of the I-population recruited is roughly
commensurate with the fraction of E-neurons participating in the event, and 
E-firing may lead by a little bit as the event is initiated by the crossing of threshold for
a few E-cells. 
The tight relationship between our dashed purple
and green curves in Fig 5 captured well  this phenomenon.

\section{Control of frequencies and spectrum}

For periodic signals, Fourier coefficients provide the right mathematical tool for extracting the main frequencies. To capture the pseudo-periodicity of gamma rhythms, neuroscientists have used
the following computational tool \cite{Hen-2005,Cha-2018}, which we review below;
it is useful but perhaps not standard for mathematicians. 
The idea is to fix a time interval of suitable length $T$ and 
compute Fourier coefficients on $[0,T]$ as if the signal was periodic with period $T$. 
This is repeated with many time shifts, i.e., we sample the signal on 
$[t, t+T]$ for $t=t_0, t_0+dt, t_0+2dt, \cdots$ for small $dt$, 
and the computed Fourier coefficients are averaged over all of these samples.
To capture gamma-band frequencies,
the time interval $T$ is usually chosen to be between $200$ and $500$ ms: too short of
an interval will fail to capture the relevant frequencies, and too long of an interval is ineffective
since the signal is not truly periodic.

\begin{figure}[t!]
 \label{fig:FTMoy}
\includegraphics[height=4cm,width=10cm]{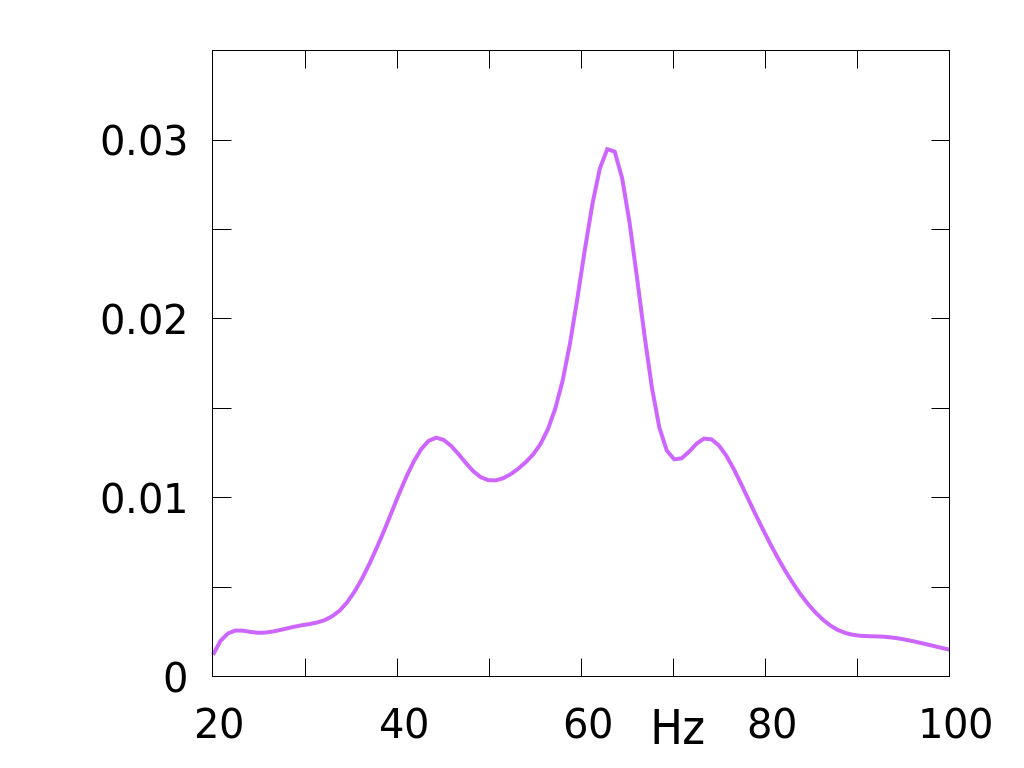}\\
\includegraphics[height=4cm,width=10cm]{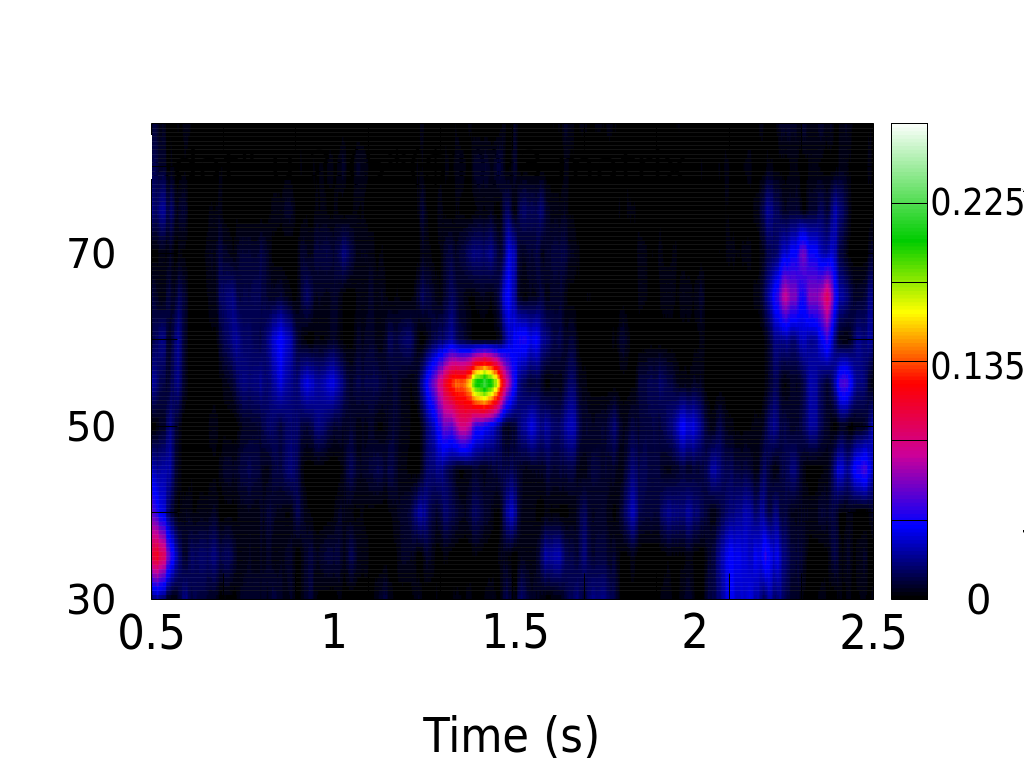}\\
\caption{Panel A.  Power spectral density (PSD) of a signal generated from Eq. (1).
The plot shows the mean values of the squares of the Fourier coefficients of  $v(t)$ computed along the time interval $[0.5s, 2.5s]$ with periods of $200 ms$. That is, Fourier coefficients are computed over time intervals $[t_0+idt,t_0+idt+T]$ with $i\in \{0,20000\}$, $t_0=500ms$, $T=200 ms$ and $dt=0.1$. The original signal comes from simulation of the stochastic system \eqref{eq:a-y}. The main frequency is around 65 Hz, which belongs to the $\gamma$-band frequency.
Panel B. This panel illustrates the  power of the same signal in a moving frame. 
Fourier coefficients are computed as above, but instead of averaging over $i$, the squares of
the coefficients are represented with specific colors and shown as a function of time.}
\end{figure}

The usual formula for the Fourier coefficients of a function $f$ is given by:
\begin{equation}
\label{eq:TF-k}
\hat{f}(k)=\frac{1}{T}\int_0^Tf(t)e^{-2i\pi k \frac{t}{T}}dt,
\end{equation}
 and the power concentrated at frequency $k$ is defined to be $|\hat f(k)|^2$.
With Neuroscience applications in mind, one might want to deal with frequencies in hertz.
In this case, a frequency $k$ on a $[0,T]$ period gives a frequency of $n=\frac{k}{T}$ in hertz. One can then adopt the following notation:
\begin{equation}
\label{eq:TF-n-1}
\hat{f}(nT)\ = \ \frac{1}{T}\int_0^Tf(t)e^{-2i\pi nT \frac{t}{T}}dt\  = \ \frac{1}{T}\int_0^Tf(t)e^{-2i\pi n t}dt.
\end{equation}
Since we deal with discrete signals with time step $\Delta t$, the computation turns into
\begin{eqnarray*}
\label{eq:TF-k-d}
\hat{f}(k) & =& \frac{1}{N}\sum_{j=0}^{N-1}f(j\Delta t)e^{-2i\pi k \frac{j\Delta t}{T}} , \quad
k \in \{0,1,...,N-1\}, \, N=\frac{T}{\Delta t}\\
\mbox{or} \qquad 
\hat{f}(nT) & = & \frac{1}{N}\sum_{j=0}^{N-1}f(j\Delta t)e^{-2i\pi n j\Delta t}, \quad 
n \in \{0,1/T,...,(N-1)/T\}\ .
\end{eqnarray*}
%

The power spectral density (PSD) plot for a signal generated by the system Eq. \eqref{eq:a-y} is 
shown in Fig 6A.
Here we let $f(t)= v(t)$ 
The graph shows the mean value of $|\hat f(k)|^2$ as function of $k$ averaged over 2 seconds. 
The broad-band nature of the signal is quite evident. In Panel B, we show, 
for the same signal as in Panel A, the concentration of power as a function of time instead 
of averaging over the samples
taken at different time points. This plot shows the wandering nature of
the frequencies and power of the signal from Eq. \eqref{eq:a-y}. It captures well the phenomenon
depicted in \cite{Xin-2012} and \cite{Cha-2018}.

\section{Summary and Discussion}

Our larger aims are (i) to help make biology quantitative and (ii) to bring biological 
ideas to mathematics, to dynamical systems  in particular. 

With regard to (i), our challenge was to simulate gamma-band activity. We went beyond
that to simulate the E and I-conductances of typical neurons, which according to the 
REI mechanism is responsible for producing gamma rhythms in a driven 
local population. We  proposed a two dimensional model inspired by the FHN and Leslie-Gower 
models, providing both theoretical and numerical analyses that give a panorama of the 
expected dynamical behavior in a range of parameters. We then introduced some 
stochasticity in the model that allows variation in amplitude, phase and frequencies, 
in a way that is consistent with real electric signals of the brain. Because solutions of
our stochastic model are realistic, they can be used in future theoretical studies to simulate
input currents to neurons under drive. More generally, the techniques demonstrated
are generalizable to simulations of other rhythmic activity in biological systems. 

With regard to (ii), we seek to contribute to the repertoire of dynamical systems theory by
bringing to the community's attention important phenomena from biology. 
Dynamical systems theory is well positioned to study periodic or quasi-periodic behavior (as in
KAM theory); there is also a fairly well developed theory of chaotic systems. Rhythms of the brain
are produced entirely naturally, and they are neither periodic nor chaotic, but somewhere in between. 
They occur as a result of a push-and-pull mechanism that is typical of the competition 
and balance between opposing groups of agents in biology. Dynamical systems 
will be enriched by a theory of 
 biological rhythms, and, more generally, by theories of interplay between opposing forces 
 from constituent subsystems.

%

\section{Appendix} 
 In Sect. 2.2.1, we mentioned the presence of some delicate behavior 
near the fold point $C$ (see also Fig 2A). Here we give further details on the analysis near this point. As discussed in section $2$, for $\epsilon$ very small, and the range of parameters under consideration,  the slow-fast analysis allows to locate the limit cycle: trajectories follow the attractive parts of the critical manifolds, exit the attractive cubic part of the critical manifold  through the jump point, and follow the fast fibers outside of the critical manifold. This gives the main picture. However, the point $C$ at which the trajectory exits the part of the critical manifold $u=0$ is significantly below the point $(0,-Ka_1a_2)$ at which the critical manifold becomes unstable. The trajectory follows the unstable part of the critical manifold along a significant distance, which is referred to the so called canard phenomenon. The aim of this appendix is to detail quite simple  computations estimating this distance and providing the coordinates of $C$. The final result is given by the theorem \ref{th:eps} at the end of the appendix. To facilitate the reading, we prepare the theorem with a few propositions.   In the first proposition, we rewrite system \eqref{eq:a-y} around the fold point $(0,-Ka_1a_2)$.

\begin{proposition}
Around the fold point $C=(0,-Ka_1a_2)$, system \eqref{eq:a-y} rewrites:
\begin{equation}
\label{eq:a-chv}
\left\{\begin{array}{rcl}
\dot{x}&=&x(Kx(a_1+a_2)+x^2-y)\\ 
\dot{y}&=&\epsilon\big(-Ka_1a_2(c+Ka_1a_2)-bKa_1a_2x+(c+2Ka_1a_2)y+bxy-y^2\big)\\
\end{array}
\right.
\end{equation}
\end{proposition}
\begin{proof}
The result follows from the change of variables
\[u=x,\, v=-Ka_1a_2+y.\]

\end{proof}
Next, we apply the change of variables:
\begin{equation}
\label{eq:chv2}
x=rx_2,\, y=r^2y_2,\, \epsilon=r^3.
\end{equation}
The following proposition holds:
\begin{proposition}
After change of variables  equation \eqref{eq:a-chv} writes
\begin{equation}
\label{eq:a-chv2}
\left\{\begin{array}{rcl}
\dot{x}_2&=&K(a_1+a_2)x_2^2+rx^3_2-ry_2\\ 
\dot{y}_2&=&-Ka_1a_2(c+Ka_1a_2)-bKa_1a_2rx_2+(c+2Ka_1a_2)r^2y_2+br^3x_2y_2-r^4y_2^4\\

\end{array}
\right.
\end{equation}
\end{proposition}
\begin{proof}
The result follows from the change variables \eqref{eq:chv2} and change of time $\tau=rt$. 
\end{proof}

Next, consider equation \cref{eq:a-chv2} with $r=0$, $i.e.$
\begin{equation}
\label{eq:K2r=0}
\left\{\begin{array}{rcl}
\dot{x}_2&=&K(a_1+a_2)x_2^2\\ 
\dot{y}_2&=&-Ka_1a_2(c+Ka_1a_2)\\
\end{array}
\right.
\end{equation}

 The following result, which follows from explicit computation, holds:

\begin{proposition}
 The solution of system \eqref{eq:K2r=0} is:
\begin{equation}
\label{eq:sol:K2r=0}
\begin{array}{rcl}
x_2(t)&=&\frac{1}{x_2^{-1}(0)-K(a_1+a_2)t}\\
y_2(t)&=&y_2(0)-Ka_1a_2(c+Ka_1a_2)t\\
\end{array}
\end{equation}
i.e.
\begin{equation*}
x_2(t)=\frac{1}{\frac{1}{x_2(0)}+\frac{a_1+a_2}{a_1a_2(c+Ka_1a_2)}(y_2(t)-y_2(0))}
\end{equation*}
or
\begin{equation*}
y_2(t)=y_2(0)-\frac{a_1a_2(c+Ka_1a_2)}{a_1+a_2}(\frac{1}{x_2(0)}-\frac{1}{x_2(t)}
\end{equation*}
It follows that orbits solutions of \cref{eq:sol:K2r=0} have the following properties:
\begin{enumerate}
\item Every orbit has a horizontal asymptote given by 
\[ y_2 = y_2(0)-\frac{a_1a_2(c+Ka_1a_2)}{a_1+a_2}\frac{1}{x_2(0)}\]
\item  Every orbit has a  vertical asymptote $x_2= 0^{+}$.

\end{enumerate}

\end{proposition}
Assume that a trajectory crosses the line $y=0$ $(v=-Ka_1a_2)$ at 
\[x(0)=k\epsilon +o(\epsilon)\]
Let 
\[\bar{x}=\frac{a_2}{2},\]
and let $\bar{y}$ the ordinate at which the trajectory crosses the line $x=\bar{x}$.
The following theorem holds.

\begin{theorem}
\label{th:eps}
Assume that the trajectory crosses the line $y=0$ at 
\[x(0)=k\epsilon +o(\epsilon)\]
then
\[\lim_{\epsilon \rightarrow 0}\bar{y}=-\frac{a_1a_2(c+Ka_1a_2)}{k(a_1+a_2)}\]
\end{theorem}
\begin{proof}
We have 
\begin{equation*}
\begin{array}{rcl}
\bar{y}&=&\epsilon^{\frac{2}{3}}\bar{y}_2\\[0.2cm]
&=&\epsilon^{\frac{2}{3}}\big(\frac{a_1a_2(c+Ka_1a_2)}{a_1+a_2}(\frac{1}{\bar{x}_2}-\frac{1}{x_2(0)})+\epsilon^{\frac{1}{3}}O(1) \big)\\[0.2cm]
&=&\epsilon^{\frac{2}{3}}\big(\frac{a_1a_2(c+Ka_1a_2)}{a_1+a_2}(\frac{\epsilon^{\frac{1}{3}}}{\bar{x}}-\frac{\epsilon^{\frac{1}{3}}}{k\epsilon +o(\epsilon)})+\epsilon^{\frac{1}{3}}O(1) \big)\\[0.2cm]
&=&\frac{a_1a_2(c+Ka_1a_2)}{a_1+a_2}(\frac{\epsilon}{\bar{x}}-\frac{1}{k +h(\epsilon)})+\epsilon O(1)  \mbox{ where } \lim_{\epsilon\rightarrow 0} h(\epsilon)=0\\
\end{array}
\end{equation*}
which gives the result.
\end{proof}
\section*{Acknowledgments} BA would like to thank R\'egion Normandie France and ERDF (European Regional Development Fund) XTERM, and CNRS IEA for funding. Part of the research of LSY was funded by NSF Grant
1901009.

\bibliographystyle{bmc-mathphys} 
\bibliography{references}      




\end{document}